\documentclass[reqno,11pt,oneside]{amsart}
\usepackage{amsfonts}
\usepackage{amsmath}
\usepackage{amsthm, amscd}
\usepackage{mathtools}
\usepackage{amssymb}
\usepackage{mathrsfs}
\usepackage{graphicx}
\usepackage{mathabx}
\usepackage{bbm}
\usepackage{overpic}

\usepackage[alphabetic]{amsrefs}
\usepackage{etex}
\usepackage[hidelinks]{hyperref}

\usepackage{xypic}
\usepackage{geometry}
\usepackage{fancyhdr}
\usepackage{color}
\usepackage{overpic}
\usepackage{tikz-cd}

\allowdisplaybreaks

\theoremstyle{plain}\newtheorem{thm}{Theorem}[section]
\theoremstyle{plain}\newtheorem{cor}[thm]{Corollary}
\theoremstyle{plain}\newtheorem{prop}[thm]{Proposition}
\theoremstyle{plain}\newtheorem{lem}[thm]{Lemma}
\theoremstyle{plain}
\theoremstyle{plain}
\theoremstyle{plain}
\theoremstyle{plain}
\theoremstyle{plain}
\theoremstyle{plain}

\theoremstyle{definition}
\newtheorem{defn}[thm]{Definition}

\newtheorem{exmp}[thm]{Example}

\theoremstyle{remark}
\newtheorem{rem}[thm]{Remark}

\makeatletter
\newtheorem*{rep@theorem}{\rep@title}
\newcommand{\newreptheorem}[2]{%
	\newenvironment{rep#1}[1]{%
		\def\rep@title{#2 \ref{##1}}%
		\begin{rep@theorem}}%
		{\end{rep@theorem}}}
\makeatother
\newreptheorem{theorem}{Theorem}

\numberwithin{equation}{section}

\DeclareMathOperator{\rank}{rank}
\DeclareMathOperator{\II}{I}

\DeclareMathOperator{\SHI}{SHI}
\DeclareMathOperator{\AHI}{AHI}

\DeclareMathOperator{\SiKh}{\Sigma Kh}
\DeclareMathOperator{\AKh}{AKh}

\DeclareMathOperator{\CKh}{CKh}

\DeclareMathOperator{\spann}{span}

\DeclareMathOperator{\id}{id}

\DeclareMathOperator{\SiHI}{\Sigma HI}

\newcommand{\bC}{\mathbb{C}}

\newcommand{\bZ}{\mathbb{Z}}

\newcommand{\bfv}{\mathbf{v}}

\newcommand{\be}{\beta}
\newcommand{\ga}{\gamma}

\newcommand{\ot}{\otimes}

\newcommand{\bpf}{\begin{proof}}
\newcommand{\epf}{\end{proof}}
\newcommand{\bthm}{\begin{thm}}
\newcommand{\ethm}{\end{thm}}
\newcommand{\bprop}{\begin{prop}}
\newcommand{\eprop}{\end{prop}}
\newcommand{\bcor}{\begin{cor}}
\newcommand{\ecor}{\end{cor}}
\newcommand{\blem}{\begin{lem}}
\newcommand{\elem}{\end{lem}}
\newcommand{\bdefn}{\begin{defn}}
\newcommand{\edefn}{\end{defn}}
\newcommand{\bexmp}{\begin{exmp}}
\newcommand{\eexmp}{\end{exmp}}
\newcommand{\brem}{\begin{rem}}
\newcommand{\erem}{\end{rem}}

\newcommand{\bdia}{\begin{displaymath}\xymatrix}
\newcommand{\edia}{\end{displaymath}}
\newcommand{\beq}{\begin{equation*}\begin{aligned}}
\newcommand{\eeq}{\end{aligned}\end{equation*}}

\newcommand{\mfo}{\mathfrak{o}}

\usepackage{extarrows}

\author{Zhenkun Li}
\address{Department of Mathematics, Stanford University, California 94305, USA}
\email{zhenkun@stanford.edu}
\author{Yi Xie}
\address{Beijing International Center for Mathematical Research, Peking University, Beijing 100871, China}
\email{yixie@pku.edu.cn}
\author{Boyu Zhang}
\address{Mathematics Department, University of Maryland at College Park, Maryland 20742, USA}
\email{bzh@umd.edu}

\title{A deformation of Asaeda--Przytycki--Sikora homology}

\begin{document}

\begin{abstract}
	We define a $1$-parameter family of homology invariants for links in thickened oriented surfaces. It recovers the homology invariant of Asaeda--Przytycki--Sikora \cite{APS} and the invariant defined by Winkeler \cite{winkeler2021khovanov}.  The new invariant can be regarded as a deformation of Asaeda--Przytycki--Sikora homology; it is not a Lee--type deformation as the deformation is only non-trivial when the surface is not simply connected. Our construction is motivated by computations in singular instanton Floer homology. We also prove a detection property for the new invariant, which is a stronger result than the main theorem of \cite{li2022instanton}.
\end{abstract}

\maketitle

\section{Introduction}
Khovanov homology \cite{Kh-Jones} is a link invariant that assigns a bi-graded homology group to every oriented link in $\mathbb{R}^3$. Asaeda--Przytycki--Sikora \cite{APS} introduced a generalization of Khovanov homology for links in $(-1,1)$--bundles over surfaces, where the bundles are required to be oriented as $3$--manifolds.  Such $(-1,1)$--bundles are called \emph{thickened surfaces}. When the surface is an annulus, Asaeda--Przytycki--Sikora homology is also called \emph{annular Khovanov homology}. Khovanov homology and Asaeda--Przytycki--Sikora homology have been essential tools for the study of knots and links for decades. More recently, Winkeler \cite{winkeler2021khovanov} introduced another variation of Khovanov homology for links in thickened multi-punctured disks, which is different from the invariant of Asaeda--Przytycki--Sikora.

Suppose $\Sigma$ is an oriented surface. In this paper, we define a one-parameter family of homology invariants for oriented  links in $(-1,1)\times \Sigma$. As bi-graded modules, the new invariant recovers both Asaeda--Przytycki--Sikora homology and the invariant of Winkeler, and it can be interpreted as a one-parameter deformation of Asaeda--Przytycki--Sikora homology. The deformation is not a Lee--type deformation as it is only non-trivial when the surface has a non-trivial fundamental group. 
The construction is motivated by computations from singular instanton Floer homology. We also use instanton Floer theory to prove a detection result for the deformed Asaeda--Przytycki--Sikora homology, which gives a stronger rank estimate than the main theorem of \cite{li2022instanton}.

The paper is organized as follows. Section \ref{sec_notation} introduces some notation and conventions. Sections \ref{sec_band} and \ref{sec_commutativity} define the differential map and proves that $d^2=0$. Section \ref{sec_Khovanov} defines the homology invariant and proves the invariance under Reidemeister moves. Section \ref{sec_motivation} explains the motivation from instanton Floer homology and prove the aforementioned detection result in Theorem \ref{thm_detection}.

\section{Notation}
\label{sec_notation}
Throughout this paper, we use $R$ to denote a fixed commutative ring with unit.  
We use $\Sigma$ to denote an oriented surface, possibly with boundary and possibly non-compact. 

For every embedded closed 1-manifold $c\subset \Sigma$, we assign an $R$--module $V(c)$ to $c$ as follows:
\begin{enumerate}
	\item  If $\ga$ is a contractible simple closed curve on $\Sigma$, define $V(\ga)$ to be the free $R$--module generated by $\bfv(\ga)_+$ and $\bfv(\ga)_-$, where $\bfv(\ga)_+$ and $\bfv(\ga)_-$ are formal generators associated with $\ga$. 
	\item If $\ga$ is a non-contractible simple closed curve, let $\mathfrak{o}$, $\mathfrak{o}'$ be the two orientations of $\ga$. Define $V(\ga)$ to be the free module generated by $\bfv(\ga)_{\mfo}$ and $\bfv(\ga)_{\mfo'}$, where $\bfv(\ga)_{\mfo}$ and $\bfv(\ga)_{\mfo'}$ are formal generators.
	\item  In general, suppose the connected components of $c$ are $\ga_1,\dots,\ga_k$, define $V(c)$ to be 
	$ \ot_{i=1}^k V(\ga_i).$
\end{enumerate}

When the choice of $\Sigma$ needs to be emphasized, we will write $V(c)$ as $V^\Sigma(c)$, and write  $\bfv(\ga)_{\mfo}$, $\bfv(\ga)_\pm$ respectively as $\bfv^\Sigma(\ga)_{\mfo}$, $\bfv^\Sigma(\ga)_\pm$.

If $\mfo$ is an orientation of a curve $\ga$, we use $\ga_\mfo$ to denote the corresponding oriented curve.

\section{Band surgery homomorphisms}
\label{sec_band}

\begin{figure}
	\begin{overpic}[width=0.8\textwidth]{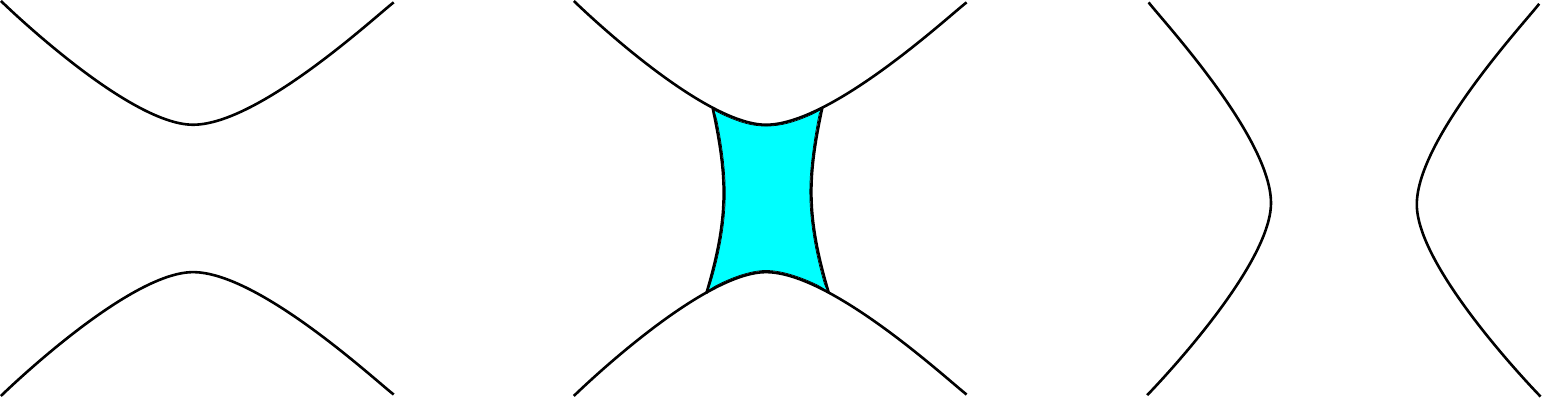}
		\put(2,-5){the manifold $c$}
		\put(42,-5){the band $b$}
		\put(78,-5){the manifold $c_b$}
	\end{overpic}
	\vspace{\baselineskip}
	\caption{Band surgery}\label{fig:band_surgery}
\end{figure}

Suppose $c$ is an embedded closed $1$--manifold on $\Sigma$, suppose $b$ is an embedded disk on $\Sigma$ such that the interior of $b$ is disjoint from $c$ and the boundary of $b$ intersects $c$ at two arcs (see Figure \ref{fig:band_surgery}). The surgery of $c$ along $b$ yields another embedded closed $1$--manifold on $\Sigma$, which we denote by $c_b$.  We will call the disk $b$ a \emph{band} that is \emph{attached} to $c$.

For later reference, we record the following two elementary lemmas.
\begin{lem}
	\label{lem_band_surgery_three_cases}
	The change from $c$ to $c_b$ has three possibilities:
	\begin{enumerate}
		\item two circle components of $c$ are merged to one circle,
		\item one circle component of $c$ is split to two circles,
		\item one circle component of $c$ is modified by the surgery to another circle.
	\end{enumerate}
\end{lem}

\begin{proof}
	Since $\partial b\cap c$ contains two arcs, at most two components of $c$ are affect by the surgery. If $\partial b\cap c$ are on two different components of $c$, then the surgery merges these two components into one circle. If $\partial b \cap c$ are on one component of $c$, then the boundary orientation of $b$ defines an orientation on both components of $\partial b\cap c$, so we have two oriented arcs embedded in one component $\ga$ of $c$.  If these two arcs induce the same orientation on $\ga$, then the surgery splits one component of $c$ to two circles. If these two arcs induce opposite orientations on $\ga$, then the surgery changes this component to another circle.
\end{proof}

Recall that if $\mfo$ is an orientation of a curve $\ga$, we use $\ga_\mfo$ to denote the corresponding oriented curve. 

\begin{lem}
	\label{lem_orientation_reverse_not_isotopy}
	Suppose $\ga$ is a simple closed curve on a connected surface $\Sigma$, and assume $\Sigma$ is not diffeomorphic to $S^2$. Suppose $\mfo$, $\mfo'$ are the two orientations of $\ga$. Then $\ga_\mfo$ and $\ga_{\mfo'}$ are not isotopic on $\Sigma$.
\end{lem}
\begin{proof}
	
	If $\ga$ is non-separating, there exists an oriented simple closed curve $\beta$ such that the algebraic intersection number of $\be$ and $\ga$ is non-zero. Since isotopies preserve the sign of algebraic intersection numbers, the desired result follows. 
	
	If $\ga$ is separating and $\partial \Sigma\neq \emptyset$, then every orientation of $\ga$ defines an ordering of the two components of $\Sigma\backslash \ga$, which defines an ordered partition of the components of $\partial \Sigma$. Since every isotopy of $\ga$ on $\Sigma$ can be extended to an isotopy of $\Sigma$ fixing the boundary, the desired result is proved.
	
	If $\ga$ is separating and $\Sigma$ is closed, then every orientation of $\ga$ defines an ordering of the two components of $\Sigma\backslash \ga$. Suppose $\Sigma_1$ and $\Sigma_2$ are the two components of $\Sigma\backslash \ga$ ordered by an orientation $\mfo$ of $\ga$. Since $\Sigma$ is not a sphere, the images of $H_1(\Sigma_1;\bZ)$ and $H_1(\Sigma_2;\bZ)$ are distinct in $H_1(\Sigma;\bZ)$. The images of $H_1(\Sigma_1;\bZ)$ and $H_1(\Sigma_2;\bZ)$ are invariant under isotopies of $\ga_\mfo$, so the desired result is proved.
\end{proof}

Take an arbitrary element $\lambda\in R$, we define a homomorphism 
$$T_\lambda(b):V(c)\to V(c_b)$$ associated with the band surgery along $b$. When the choice of $\Sigma$ needs to be emphasized, we will write $T_\lambda(b)$ as $T^\Sigma_\lambda(b)$.

We first assume that the intersection of $\partial b$ with every component of $c$ is non-empty. The general case will be discussed later. By Lemma \ref{lem_band_surgery_three_cases}, if the intersection of $\partial b$ with every component of $c$ is non empty, then there are three cases:

{\bf Case 1}: $c$ has two components $\ga_1$ and $\ga_2$ and they are merged into one circle $\ga=c_b$ after the surgery. In this case, we define 
$T_\lambda(b): V(\ga_1)\ot V(\ga_2)\to V(\ga)$
as follows:
\begin{enumerate}
	\item If both $\ga_1$ and $\ga_2$ are contractible circles, then $\ga$ is also contractible, and we define $T_\lambda(b)$ by
	\begin{align*}
		\mathbf{v}(\ga_1)_+ \otimes \mathbf{v}(\ga_2)_+ &\mapsto \mathbf{v}(\ga)_+,  &\mathbf{v}(\ga_1)_+& \otimes \mathbf{v}(\ga_2)_- \mapsto \mathbf{v}(\ga)_- ,\\
		\mathbf{v}(\ga_1)_- \otimes \mathbf{v}(\ga_2)_+ &\mapsto \mathbf{v}(\ga)_-,  &\mathbf{v}(\ga_1)_-& \otimes \mathbf{v}(\ga_2)_- \mapsto 0.
	\end{align*}
	\item If $\gamma_1$  is contractible and $\gamma_2$ is non-contractible, then $\ga_2$ is isotopic to $\ga$. The existence of non-contractible curves on $\Sigma$ implies that $\Sigma$ is not diffeomorphic to $S^2$. By Lemma \ref{lem_orientation_reverse_not_isotopy}, the orientations of $\ga_2$ are canonically identified with the orientations of $\ga$ via an isotopy. This identification defines a canonical isomorphism from $V(\ga_2)$ to $V(\ga)$, which we denote by $\iota$. 
	In this case, the homomorphism $T_\lambda(b)$ is defined by
	\[
	\mathbf{v}(\gamma_1)_+\otimes x \mapsto \iota(x),\quad  \bfv(\ga_1)_-\otimes x \mapsto 0
	\]
	for all $x\in V(\ga_2)$. 
	\item If $\gamma_1$  is non-contractible and $\gamma_2$ is contractible, define $T_\lambda(b)$ by requiring the map to be symmetric with respect to $\ga_1$ and $\ga_2$ and deducing to case (2) above.
	\item If $\gamma_1$ and $\gamma_2$ are both non-contractible and $\gamma_3$ is contractible, then $\ga_1$ and $\ga_2$ must be isotopic. 
	By Lemma \ref{lem_orientation_reverse_not_isotopy}, the orientations of $\ga_1$ and $\ga_2$ are canonically identified by the isotopy. Let $\mfo$, $\mfo'$ be the two orientations of $\ga_1$, and  use the same notation to denote the corresponding orientations of $\ga_2$. 
	The map $T_\lambda(b)$ is then defined by
	\begin{align*}
		\mathbf{v}(\gamma_1)_\mfo\otimes \mathbf{v}(\gamma_2)_\mfo &\mapsto 0,   &\mathbf{v}(\gamma_1)_{\mfo'}& \otimes \mathbf{v}(\gamma_2)_{\mfo'}\mapsto 0,\\
		\mathbf{v}(\gamma_1)_{\mfo}\otimes \mathbf{v}(\gamma_2)_{\mfo'} &\mapsto \mathbf{v}(\gamma)_-, &\mathbf{v}(\gamma_1)_{\mfo'}&\otimes \mathbf{v}(\gamma_2)_{\mfo} \mapsto \mathbf{v}(\gamma)_-.
	\end{align*}
	\item If all of $\ga_1$, $\ga_2$, and $\ga$ are non-contractible, let $N$ be the regular neighborhood of $b\cup \ga_1 \cup \ga_2$. Then $N$ is a sphere with three disks removed, and the three boundary components of $N$ are parallel to $\ga_1$, $\ga_2$, $\ga$.  Since $N\subset \Sigma$ is oriented, the boudary orientation of $N$ defines an orientation on each of $\ga_1, \ga_2, \ga$, and we denote these orientations by $\mfo_1,\mfo_2,\mfo$ respectively.  Denote their opposite orientations by $\mfo_1',\mfo_2',\mfo'$. Then $T_\lambda(b)$ is defined by
	\begin{align*}
		\mathbf{v}(\gamma_1)_{\mfo_1'}\otimes \mathbf{v}(\gamma_2)_{\mfo_2'} &\mapsto  \lambda\cdot \mathbf{v}(\gamma)_{\mfo},   &\mathbf{v}(\gamma_1)_{\mfo_1'}& \otimes \mathbf{v}(\gamma_2)_{\mfo_2}\mapsto 0,\\
		\mathbf{v}(\gamma_1)_{\mfo_1}\otimes \mathbf{v}(\gamma_2)_{\mfo_2'} &\mapsto 0, &\mathbf{v}(\gamma_1)_{\mfo_1}&\otimes \mathbf{v}(\gamma_2)_{\mfo_2} \mapsto 0.
	\end{align*}
\end{enumerate}

{\bf Case 2}: $c$ contains one component $\ga$ and $c_b$ has two components $\ga_1$ and $\ga_2$. In this case, we define  $T_\lambda(b): V(\ga)\to V(\ga_1)\ot V(\ga_2)$ as follows:
\begin{enumerate}
	\item If $\gamma_1$ and $\gamma_2$ are both contractible circles, then $\ga$ is also contractible, and we define $T_\lambda(b)$ by
	\begin{align*}
		\mathbf{v}(\gamma)_+  &\mapsto \mathbf{v}(\gamma_1)_+\otimes \mathbf{v}(\gamma_2)_- + \mathbf{v}(\gamma_1)_-\otimes \mathbf{v}(\gamma_2)_+, \\
		\mathbf{v}(\gamma)_-  &\mapsto \mathbf{v}(\gamma_1)_-\otimes \mathbf{v}(\gamma_2)_-.
	\end{align*}
	
	\item If one of $\{\ga_1,\ga_2\}$ is contractible and the other is non-contractible, assume without loss of generality that $\ga_1$ is contractible and $\ga_2$ is non-contractible. Then $\ga$ is isotopic to $\ga_2$, and the orientations of $\ga$ and $\ga_2$ are canonically identified.  Let $\mfo$, $\mfo'$ be the two orientations of $\ga$, and use the same notation to denote the corresponding orientations of $\ga_2$. Define the map $T_\lambda(b)$ by 
	$$
	\mathbf{v}(\gamma)_{\mfo} \mapsto \mathbf{v}(\gamma_1)_-\otimes \mathbf{v}(\gamma_2)_{\mfo}, \quad \mathbf{v}(\gamma)_{\mfo'} \mapsto \mathbf{v}(\gamma_1)_-\otimes \mathbf{v}(\gamma_2)_{\mfo'}.
	$$
	\item If both $\ga_1$ and $\ga_2$ are non-contractible and $\ga$ is contractible, then $\ga_1$, $\ga_2$ are isotopic to each other, and the orientations of $\ga_1$ are $\ga_2$ are canonically identified. Let $\mfo,\mfo'$ be the orientations of $\ga_1$ and use the same notation for the orientations of $\ga_2$. Define the map $T_\lambda(b)$ by
	$$
	\mathbf{v}(\gamma)_+ \mapsto \mathbf{v}(\gamma_1)_{\mfo}\otimes \mathbf{v}(\gamma_2)_{\mfo'} + \mathbf{v}(\gamma_1)_{\mfo'}\otimes \mathbf{v}(\gamma_2)_{\mfo} , \quad \mathbf{v}(\gamma)_- \mapsto 0.
	$$
	\item If all of $\ga,\ga_1,\ga_2$ are non-contractible, let $N$ be the regular neighborhood of $b\cup \ga$.  Then $N$ is a sphere with three disk removed, and the three boundary components of $N$ are parallel to $\ga_1$, $\ga_2$, $\ga$.  The boundary orientation of $N$ defines an orientation on each of $\ga_1, \ga_2, \ga$, and we denote them by $\mfo_1,\mfo_2,\mfo$ respectively.  Denote their opposite orientations by $\mfo_1',\mfo_2',\mfo'$. Define the map $T_\lambda(b)$ by
	\[
	\mathbf{v}(\gamma)_{\mfo'}\mapsto \lambda\cdot \bfv(\ga_1)_{\mfo_1}\otimes \mathbf{v}(\gamma_2)_{\mfo_2},\quad \bfv(\ga)_{\mfo}\mapsto 0.
	\]
\end{enumerate}

{\bf Case 3}: both $c$ and $c_b$ have exactly one component. In this case, define $T_\lambda(b)$ to be zero.\\

In general, suppose $c=c^{(1)}\sqcup c^{(2)}$ such that $\partial b$ is disjoint from $c^{(2)}$ and intersects every component of $c^{(1)}$, we define the band surgery homomorphism
$T_\lambda(b) : V_\lambda(c)\to V_\lambda(c_b)$
to be 
\begin{equation}
	\label{eqn_def_band_homomorphism_general_case}
	T_\lambda(b) = T_\lambda(b)|_{V(c^{(1)})}\otimes \id|_{V(c^{(2)})}.
\end{equation}

\begin{rem}
	In the above definition, the coefficient $\lambda$ only appeared in Cases 1(5) and 2(4).
\end{rem}

\section{Commutativity of band surgery homomorphisms}
\label{sec_commutativity}
The main result of this section is the following proposition.
\begin{prop}
	\label{prop_commutativity_band_hom}
	Suppose $c$ is an embedded closed $1$--manifold on $\Sigma$, and suppose $b_1$ and $b_2$ are two disjoint bands attached to $c$.  Then for all $\lambda\in R$,
	\begin{equation}
		\label{eqn_commutativity_band_hom}
		T_\lambda(b_1)\circ T_\lambda(b_2) = T_\lambda(b_2) \circ T_\lambda(b_1).
	\end{equation}
\end{prop}

\subsection{The genus-zero case}
We first establish \eqref{eqn_commutativity_band_hom} when $\Sigma$ is a sphere or a finitely punctured sphere. Our argument here is inspired by the work of Winkeler \cite{winkeler2021khovanov}.

\begin{lem}
	\label{lem_genus_zero_commutativity}
	Equation \eqref{eqn_commutativity_band_hom} holds if $\Sigma$ is a sphere or a finitely punctured sphere.
\end{lem}

\begin{proof}
	If $\Sigma$ is a sphere or a disk, then every curve is contractible, and Case (3) in Lemma \ref{lem_band_surgery_three_cases} is not possible. In this case, our definition of $T_\lambda(b)$ does not depend on $\lambda$ and it coincides with the definition of the merge and split maps in standard Khovanov theory. Therefore Equation \eqref{eqn_commutativity_band_hom} holds.
	
	When $\Sigma$ has $n\ge 2$ boundary components, we view $\Sigma$ as a disk $B$ with $n-1$ interior disks $B_1,\dots,B_{n-1}$ removed. 
	Assume the orientation of $\Sigma$ is defined so that the boundary orientation on $\partial B$ is given by the counter-clockwise orientation, and the boundary orientation on $\partial B_i$ is the clockwise orientation.
	
	Recall that when the surface $\Sigma$ needs to be emphasized, we write $V(c)$, $\bfv(\ga)_{\mfo}$, $\bfv(\ga)_\pm$, $T_\lambda(b)$ respectively as $V^\Sigma(c)$, $\bfv^\Sigma(\ga)_{\mfo}$, $\bfv^\Sigma(\ga)_\pm$, $T_\lambda^\Sigma(b)$.
	
	For each embedded closed $1$--manifold $c\subset \Sigma$, 
	define an isomorphism $\Phi: V^B(c)\to V^\Sigma(c)$ as follows. For each component $\ga$ of $c$, if $\ga$ is contractible in $\Sigma$,  define $$\Phi(\bfv^B(\ga)_\pm )=\bfv^\Sigma(\ga)_\pm.$$ If $\ga$ is non-contractible in $\Sigma$, let $\mfo$ denote the counter-clockwise orientation of $\ga$,  let $\mfo'$ denote the clockwise orientation of $\ga$, and define 
	$$\Phi(\bfv^B(\ga)_+) = \bfv^\Sigma(\ga)_\mfo,\quad \Phi(\bfv^B(\ga)_-) = \bfv^\Sigma(\ga)_{\mfo'}.$$ 
	Since $T_\lambda^B(b)$ does not depend on $\lambda$, we denote it by $T^B(b)$.
	Then 
	$$ \Phi\circ T^B(b)\circ \Phi^{-1}$$
	is a homomorphism from $V^\Sigma(c)$ to $V^\Sigma(c_b)$. 
	
	For each $i\in \{1,\dots,n-1\}$, define a grading on $V^\Sigma(c)$ as follows. If a circle $\ga$ is a contractible curve on $\Sigma$, define the degree of $\bfv^\Sigma(\ga)_\pm$ to be zero. If $\ga$ is non-contractible, for each orientation $\mfo$ of $\ga$, define the degree of $\bfv^\Sigma(\ga)_\mfo$ to be the rotation number of $\ga_\mfo$ around $B_i$. Here, our convention on the rotation number is defined so that counter-clockwise orientations always have non-negative rotation numbers. Define the grading of the tensor product of a set of generators to be the sum of the grading of each generator. 
	
	By checking all the cases in the definition of $T_\lambda(b)$, it is straightforward to verify that the map $T^\Sigma(b)$ preserves all the $n-1$ gradings defined above. Moreover, for each $i\in \{1,\dots,n-1\}$, the map $ \Phi \circ T^B(b)\circ \Phi^{-1}$ does not increase the $i^{th}$ grading. The components of $ \Phi \circ T^B(b)\circ \Phi^{-1}$ that preserve all the $n-1$ gradings is equal to the map $T^\Sigma_1(b)$, which is the map $T^\Sigma_\lambda$ when $\lambda=1$. Since  $T^B(b_1)\circ T^B(b_2) = T^B(b_2) \circ T^B(b_1)$ on $B$, we conclude that \eqref{eqn_commutativity_band_hom}  holds for $T^\Sigma_1$.
	
	To show that \eqref{eqn_commutativity_band_hom}  holds for general $\lambda$, define $T_\delta^\Sigma = T_1^\Sigma - T_0^\Sigma$. Then 
	$$
	T_\lambda^\Sigma = T_0^\Sigma + \lambda\cdot T_\delta^\Sigma.
	$$
	We define another grading on $V^\Sigma(-)$ as follows.  
	If a circle $\ga$ is a contractible curve on $\Sigma$, define the degree of $\bfv_\Sigma(\ga)_\pm$ to be zero. If $\ga$ is non-contractible, for each orientation $\mfo$ of $\ga$, define the degree of $\bfv^\Sigma(\ga)_\mfo$ to be $1$ if $\mfo$ is the counter-clockwise orientation, and define the degree of $\bfv^\Sigma(\ga)_\mfo$ to be $-1$ if $\mfo$ is the clockwise orientation. Define the grading of the tensor product of a set of generators to be the sum of the grading of each generator. 
	
	By checking all the cases in the definition of $T_\lambda^\Sigma$, it is straightforward to verify that under the above grading, the map $T_0^\Sigma$ is homogeneous with degree $0$, and $T_\delta^\Sigma$ is homogeneous with degree $-1$. Since \eqref{eqn_commutativity_band_hom} holds for $\lambda=1$, we have
	\begin{align*}
		T_0^\Sigma(b_1) \circ T_0^\Sigma(b_2) &=  T_0^\Sigma(b_2) \circ T_0^\Sigma(b_1) \\
		T_\delta^\Sigma(b_1) \circ T_0^\Sigma(b_2) + T_0^\Sigma(b_1) \circ T_\delta^\Sigma(b_2) &=  	T_\delta^\Sigma(b_2) \circ T_0^\Sigma(b_1) + T_0^\Sigma(b_2) \circ T_\delta^\Sigma(b_1) \\
		T_\delta^\Sigma(b_1) \circ T_\delta^\Sigma(b_2) &=  T_\delta^\Sigma(b_2) \circ T_\delta^\Sigma(b_1) 	
	\end{align*}
	Therefore \eqref{eqn_commutativity_band_hom} holds for all $\lambda\in R$.
\end{proof}

Lemma \ref{lem_genus_zero_commutativity} can be used to verify \eqref{eqn_commutativity_band_hom}  on surfaces with positive genera because of the following lemma.
\begin{lem}
	\label{lem_embedded_surface}
	Suppose $\Sigma$ is an oriented surface, and $\Sigma'\subset \Sigma$ is an embedded surface whose orientation is induced by $\Sigma$. Suppose the embedding of $\Sigma'$ in $\Sigma$ is $\pi_1$--injective. Suppose $c$ is an embedded closed $1$--manifold in $\Sigma'$ and $b_1,b_2$ are two disjoint bands in $\Sigma'$ attached to $c$. Then 
	$$
	T^{\Sigma'}_\lambda(b_1)\circ T^{\Sigma'}_\lambda(b_2) = T^{\Sigma'}_\lambda(b_2) \circ T^{\Sigma'}_\lambda(b_1)
	$$
	on $V_{\Sigma'}(c)$ if and only if 
	$$
	T_{\Sigma}(b_1)\circ T_{\Sigma}(b_2) = T_{\Sigma}(b_2) \circ T_{\Sigma}(b_1)
	$$
	on $V_{\Sigma}(c)$.
\end{lem}

\begin{proof}
	Since the embedding of $\Sigma'$ in $\Sigma$ is $\pi_1$--injective, there is a canonical isomorphism from $V_{\Sigma'}(c)$ to $V_{\Sigma}(c)$ for every embedded $1$--manifold $c\subset \Sigma'$ which takes the generators of $V_{\Sigma'}(c)$ to the corresponding generators of $V_{\Sigma}(c)$, and this isomorphism intertwines with $T^{\Sigma'}_\lambda$ and $T^\Sigma_\lambda$, so the lemma is proved. 
\end{proof}

\subsection{The genus-one case}
Now we prove Proposition \ref{prop_commutativity_band_hom} when $\Sigma$ is a torus or a finitely punctured torus. Let $\Sigma_0$ be a torus and suppose $\Sigma = \Sigma_0\backslash\{p_1,\dots,p_n\}$ with $n\ge 0$. Let $c,b_1,b_2$ be as in Proposition \ref{prop_commutativity_band_hom}. By the definition of $T_\lambda$, we may assume without loss of generality that every component of $c$ intersects $\partial (b_1\cup b_2)$ non-trivially.

\begin{figure}
	\begin{overpic}[width=0.8\textwidth]{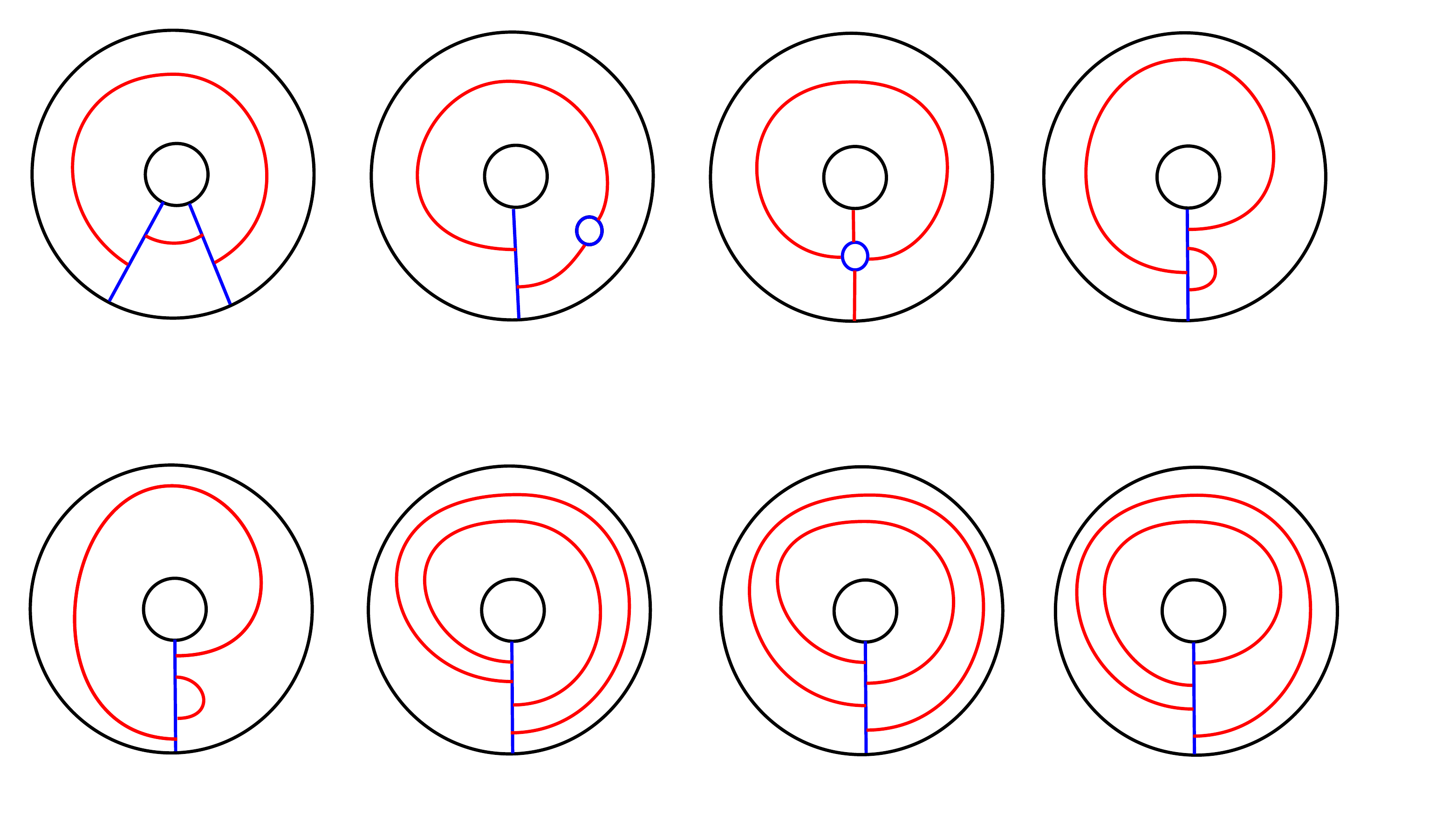}
		\put(10,0){(5)}
		\put(33,0){(6)}
		\put(58,0){(7)}
		\put(79,0){(8)}
		\put(10,30){(1)}
		\put(33,30){(2)}
		\put(58,30){(3)}
		\put(79,30){(4)}
	\end{overpic}
	\caption{All possible configurations}\label{fig:8_configurations}
\end{figure}
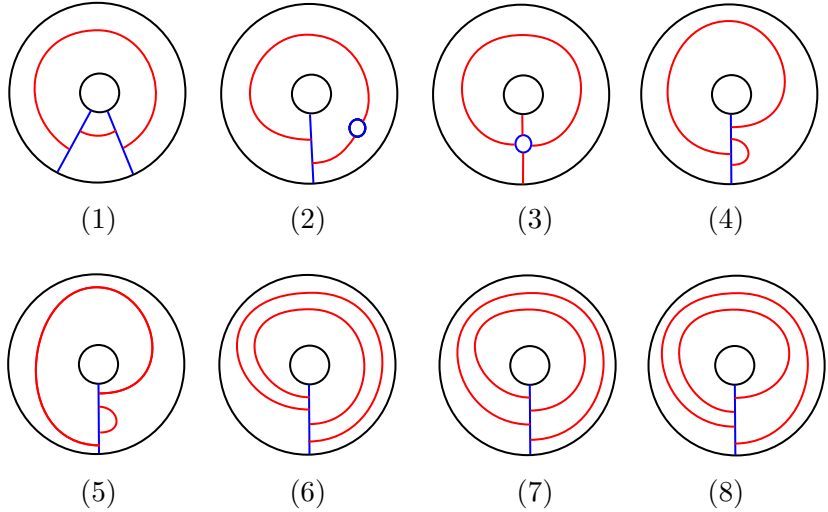

\begin{lem}
	\label{lem_eight_cases_genus_one}
	Assume every simple closed curve $\ga_0\subset \Sigma_0$ that is disjoint from $c\cup b_1\cup b_2$ is contractible in $\Sigma_0$. Then up to orientation-preserving diffeomorphisms of $\Sigma_0$, there are only 8 possible configurations of $c, b_1,b_2$ as subsets of $\Sigma_0$, which are shown in Figure \ref{fig:8_configurations}.
\end{lem}

In each case of Figure \ref{fig:8_configurations}, the torus $\Sigma_0$ is the quotient space obtained by gluing the two boundary components of the annulus. The blue curves denote the $1$--manifold $c$, and the disks $b_1$ and $b_2$ are defined to be the thickening of the red arcs.

\begin{proof}
	We discuss the following cases:
	
	If $c$ contains two circles $\ga_1$, $\ga_2$, and both of them are contractible, let $D_1, D_2\subset \Sigma$ denote the disks bounded by $\ga_1$, $\ga_2$. Then $D_1\cup D_2\cup b_1\cup b_2$ is a disk or an annulus, and hence there exists a circle $\ga_0$ in the complement of $c\cup b_1\cup b_2$  that is contractible, contradicting the assumptions.

	If $c$ contains two circles $\ga_1$, $\ga_2$, such that both $\ga_1$ and $\ga_2$ are non-contractible, then $\ga_1$ and $\ga_2$ must be parallel to each other. The complement $\Sigma_0\backslash(\ga_1\cup \ga_2)$ contains two components. If every simple closed curve in $\Sigma_0\backslash (c\cup b_1\cup b_2)$ is contractible in $\Sigma_0$, then the interior of $b_1$ and $b_2$ must be contained in different components of $\Sigma_0\backslash(\ga_1\cup \ga_2)$, and $\partial b_i$ must intersect both components of $c$ for each $i$. Therefore, up to orientation-preserving diffeomorphisms of $\Sigma_0$, the configuration is given by Case (1) of Figure \ref{fig:8_configurations}.

	If $c$ contains two circles $\ga_1$, $\ga_2$, where $\ga_1$ is contractible and $\ga_2$ is not contractible, let $D_1$ be the disk bounded by $\ga_1$. If either $b_1$ or $b_2$ is contained in $D_1$, then $D_1\cup c\cup b_1\cup b_2$ deformation retracts onto $\ga_2$, so there exists a non-contractible simple closed curve in $\Sigma_0$ that is disjoint from $D_1\cup c\cup b_1\cup b_2$, which contradicts the assumptions. Therefore, both $b_1$ and $b_2$ must be on the outside of $D_1$, so $b_1\cup D_1\cup b_2$ deformation retracts onto an arc with both end points on $\ga_2$.  The assumptions then imply that $c\cup b_1\cup b_2$ is given by Case (2) of Figure \ref{fig:8_configurations} up to orientation-preserving diffeomorphisms of $\Sigma_0$.

	If $c$ consists of one simple closed curve $\ga$  that is contractible in $\Sigma_0$, let $D$ be the disk bounded by $\ga$, then $b_1$ and $b_2$ must be the thickening of two disjoint arcs $r_1$ and $r_2$ in $\Sigma_0\backslash D$. For $i=1,2$, let $\overline{r_i}$ be the circle obtained by the union of $r_i$ with an arc in $D$. Since $r_1$ and $r_2$ are disjoint arcs, we may choose the arcs in $D$ so that $\overline{r_1}$  and $\overline{r_2}$ are either disjoint or intersect transversely at one point. The assumptions then imply that $\overline{r_1}$  and $\overline{r_2}$ must intersect transversely at one point. Hence the configuration is given by Case (3) of Figure \ref{fig:8_configurations} up to orientation-preserving diffeomorphisms of $\Sigma_0$.
	
	If $c$ consists of one non-contractible simple closed curve, then the possible configurations are given by Cases (4)-(8) of Figure \ref{fig:8_configurations}.
\end{proof}

\begin{lem}
	\label{lem_genus_one_commutativity}
	Equation \eqref{eqn_commutativity_band_hom} holds if $\Sigma$ is a torus or a finitely punctured torus.
\end{lem}

\begin{proof}
	If there exists a non-contractible simple closed curve $\ga_0\subset \Sigma_0$ that is disjoint from $c\cup b_1\cup b_2$, we may cut open $\Sigma_0$ along $\ga_0$, and the desired result follows from Lemma \ref{lem_genus_zero_commutativity} and Lemma 
	\ref{lem_embedded_surface}.
	Therefore, by Lemma \ref{lem_eight_cases_genus_one}, we only need to consider the 8 cases given by Figure \ref{fig:8_configurations}.
	
	In cases (2), (4), (5), (6), (7), (8), both sides of \eqref{eqn_commutativity_band_hom} are zero because Case (3) of Lemma \ref{lem_band_surgery_three_cases} appears on both sides of the equations. 
	
	For Cases (1) and (3), the complement $\Sigma\backslash (c\cup b_1\cup b_2)$ has two connected components. Therefore, by Lemma \ref{lem_embedded_surface} again, we only need to consider the cases when there is at most one puncture on each component. 
	
	Recall that $n$ denotes the number of punctures on $\Sigma_0$.
	For Case (1) with $n=0$ or $2$, and for Case (3), there is an orientation-preserving diffeomorphism of $\Sigma_0$ that preserves $c$ and $\Sigma$, is orientation-preserving on $c$, and switches $b_1$ and $b_2$. Therefore \eqref{eqn_commutativity_band_hom} holds.
	
	For Case (1) with $n=1$, it is straightforward to verify that both sides of \eqref{eqn_commutativity_band_hom} are zero.
\end{proof}

\subsection{Proof of Proposition \ref{prop_commutativity_band_hom}}
Now we prove Proposition \ref{prop_commutativity_band_hom} for the general case. 

Without loss of generality, we may assume that every component of $c$ intersects $b_1$ and $b_2$ non-trivially, and that $c\cup b_1\cup b_2$ is connected. 

In this case, $c\cup b_1\cup b_2$ is homotopy equivalent to the wedge sum of three circles.  Therefore its Euler characteristic is $-2$. 

Let $N$ be a closed regular neighborhood of $c\cup b_1\cup b_2$ in $\Sigma$. Let $\Sigma'$ be obtained from $N$ as follows: for each component $\ga$ of $\partial N$, if $\ga$ is contractible in $\Sigma$ but not contractible in $N$, then $\ga$ bounds a disk $D_\ga$ in $\Sigma$ such that $D\cap N=\ga$. Define $\Sigma'$ to be the union of $N$ and all disks $D_\ga$ as above. Then the embedding of $\Sigma'$ in $\Sigma$ is $\pi_1$--injective. Since $\chi(N)=-2$, the genus of $\Sigma'$ is $0$ or $1$. Therefore by the previous results, \eqref{eqn_commutativity_band_hom}  holds on $\Sigma'$. Hence by Lemma \ref{lem_embedded_surface}, the desired equation also holds on $\Sigma$.

\section{Khovanov homology}
\label{sec_Khovanov}

Suppose $L\subset (-1,1)\times \Sigma$ is a link. For each $\lambda$, we define a homology invariant for $L$ using the maps $T_\lambda$. 

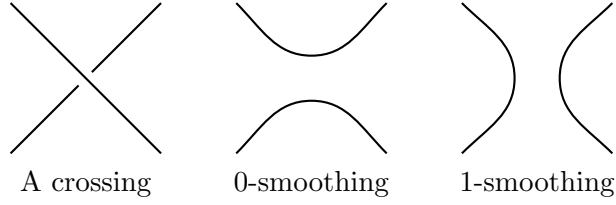
\begin{figure}
	\begin{tikzpicture}
		\draw[thick] (1,-1) to (-1,1); \draw[thick,dash pattern=on 1.3cm off 0.25cm] (1,1) to (-1,-1);  \node[below] at (0,-1.1) {A crossing};
		\draw[thick] (2,1)  to [out=315,in=180]  (3,0.3) to [out=0,in=225]   (4,1);
		\draw[thick] (2,-1)  to [out=45,in=180]  (3,-0.3) to [out=0,in=135]   (4,-1);  \node[below] at (3,-1.1) {0-smoothing};
		\draw[thick] (5,1)  to [out=315,in=90]  (5.7,0) to [out=270,in=45]    (5,-1);  \node[below] at (6,-1.1) {1-smoothing};
		\draw[thick] (7,1)  to [out=225,in=90]  (6.3,0) to [out=270,in=135]   (7,-1);
	\end{tikzpicture}
	\caption{Two types of smoothings}\label{fig_01smoothing}
\end{figure}

Suppose a link $L$ is given by a diagram $D$ on $\Sigma$ with $k$ crossings, and fix an ordering of the crossings.
For $v=(v_1,v_2,\dots, v_k)\in \{0,1\}^k$, resolving the crossings of $D$ by a sequence of $0$--smoothings and $1$--smoothings (see Figure \ref{fig_01smoothing})  by $v$ turns $D$ to an embedded closed $1$--manifold in $\Sigma$. Denote the resolved diagram by $D_v$.

Whenever $u$ is obtained from $v$ by changing one coordinate from 0 to 1, there is a band $b$ near the crossing such that $v$ is obtained from $u$ by a band surgery along $b$. Define $d_{vu}^\lambda: V(D_v) \to V(D_u)$ to be $T_\lambda(b)$. Let $e_i$ be the $i$--th standard basis vector of $\mathbb{Z}^k$.
Define
\begin{equation*}
	\CKh_{\Sigma,\lambda}(L)=\bigoplus_{v\in \{0,1\}^k} V(D_v),
\end{equation*}
and define an endomorphisms on $\CKh_\Sigma(L)$ by
\begin{equation*}
	\mathcal{D}_{\Sigma,\lambda}= \sum_i \sum_{u-v=e_i}  (-1)^{\sum_{i<j\le c}v_j}  d_{vu}.
\end{equation*}

By \eqref{eqn_commutativity_band_hom}, we have $\mathcal{D}_{\Sigma,\lambda}^2=0$.

We define a quantum grading and a homological grading on $\CKh_{\lambda,\Sigma}(L)$ as follows. For each circle $\ga$, if $\ga$ is non-contractible, define the quantum grading on $V(\ga)$ to be zero. If $\ga$ is contractible, define the quantum grading of $\bfv(\ga)_+$ to be $1$ and the quantum grading of $\bfv(\ga)_-$ to be $-1$. This grading then extends to a grading on $	\CKh_{\lambda,\Sigma}(L)$. Define the homology grading of $V(D_v)\subset \CKh_{\lambda,\Sigma}(L)$ to be the sum of coordinates in $v$.

There is also a grading on $\CKh_{\lambda,\Sigma}(L)$ over $H_1(\Sigma;\bZ)$ defined as follows. For each circle $\ga$, if $\ga$ is contractible, define the grading on $V(\ga)$ to be zero. If $\ga$ is non-contractible, for each orientation $\mfo$ of $\ga$, define the grading of $\bfv(\ga)_\mfo$ to be the fundamental class of $\ga_\mfo$. 

Following the standard convention, we use curly brackets $\{l\}$ to denote the shifting in quantum gradings by $l$ (namely, adding the quantum grading to each homogeneous element by $l$); we use the square brackets $[l]$ to denote the shifting in homology gradings by $l$.

\begin{thm}
	The homology of 
	$$\Big(\CKh_{\lambda,\Sigma}(L)[-n_-]\{n_+-2n_-\},	\mathcal{D}_{\Sigma,\lambda}\Big)$$
	as a $\bZ\oplus\bZ\oplus H_1(\Sigma;\bZ)$ graded module is independent of the diagram or the ordering of the crossings, where $n_+$ and $n_-$ denote the number of positive and negative crossings of the diagram.
\end{thm}

\begin{proof}
	The proof is identical to the proof of the invariance of the standard Khovanov homology under Reidemeister moves in  \cite{bar2002khovanov}. Besides \eqref{eqn_def_band_homomorphism_general_case} and \eqref{eqn_commutativity_band_hom}, the only properties about the band homomorphisms $T_\lambda(b)$ needed in the proof are the following:
	\begin{enumerate}
		\item If $\ga$ is a contractible circle, then $V(\ga)$ is rank $2$ with two generators $\bfv(\ga)_\pm$.
		\item Suppose the band surgery along $b$ merges two circles $\ga_1$ and $\ga_2$ to $\ga$, where $\ga_1$ is contractible. Then $\ga_2$ and $\ga$ are isotopic, and this isotopy defines a canonical isomorphism $\iota:V(\ga_2)\to V(\ga)$. Then $T_\lambda(b)(\bfv(\ga_1)_+\otimes x) = \iota(x)$ for all $x\in V(\ga_2)$.
		\item Suppose the band surgery along $b$ splits one circle $\ga$ to circles $\ga_1$ and $\ga_2$, where $\ga_1$ is contractible. Then $\ga_2$ and $\ga$ are isotopic, and this isotopy defines a canonical isomorphism $\iota:V(\ga)\to V(\ga_2)$. Then the composition map
		$$
		V(\ga) \xlongrightarrow[]{T_\lambda(b)} V(\ga_1)\otimes V(\ga_2) \xlongrightarrow[]{/\bfv(\ga_1)_+=0} \spann\{\bfv(\ga_1)_-\}\otimes V(\ga_2)
		$$
		is given by the tensor product with $\bfv(\ga_1)_-$, where the second map above is a quotient map.
	\end{enumerate}
		The only remark worth making is that there is a typo in the definition of the ``transpose'' map in Section 3.5.5 of \cite{bar2002khovanov}. The map $\Upsilon$ on the top layer should map the \emph{quotient image} of the pair $(\beta_1,\ga_1)$ to the \emph{quotient image} of the pair $(\beta_2,\ga_2)$ \emph{such that} $\ga_1 + \tau_1\beta_1 = \ga_2 + \tau_2 \beta_2$. The italicized phrases and the last equation in the previous sentence were missing in \cite{bar2002khovanov}.
\end{proof}

\begin{defn}
	We define the homology of $$\Big(\CKh_{\lambda,\Sigma}(L)[-n_-]\{n_+-2n_-\},	\mathcal{D}_{\Sigma,\lambda}\Big)$$ as a $\bZ\oplus\bZ\oplus H_1(\Sigma;\bZ)$ module to be the Khovanov invariant of $L\subset (-1,1)\times \Sigma$, and denote it by $\SiKh_\lambda(L;R)$.
\end{defn}

\begin{rem}
	When $\lambda=0$, the differential map $\mathcal{D}_{\Sigma,\lambda}$ is identical to the differential map of Asaeda--Przytycki--Sikora homology defined in \cite{APS}. When $R=\bZ$, $\lambda=1$, and $\Sigma$ is a punctured disk, the homology $\SiKh_\lambda$ recovers the invariant defined by Winkeler \cite{winkeler2021khovanov}.
\end{rem}

\section{Motivation from instanton homology}
\label{sec_motivation}
This section explains the motivation of the definition of $T_\lambda(b)$ from instanton homology. We will also prove the following detection result:
\begin{thm}
	\label{thm_detection}
	Suppose $\Sigma$ is a surface with genus zero, and $L\subset (-1,1)\times \Sigma$ is a link. Then $\rank_{\bZ/2} \SiKh_1(L;\bZ/2) \ge 2$, and the equality holds if and only if $L$ is isotopic to an embedded knot in $\Sigma$.
\end{thm}

\begin{rem}
By a spectral sequence of Winkeler \cite[Theorem 1.3]{winkeler2021khovanov}, we have 
$$
\rank_{\bZ/2} \SiKh_0(L;\bZ/2)\ge \rank_{\bZ/2} \SiKh_1(L;\bZ/2).
$$
Therefore, Theorem \ref{thm_detection} is an improvement of \cite[Theorem 1.2]{li2022instanton}.
\end{rem}

Suppose $R$ is a \emph{closed} oriented surface, and let $L$ be a link in $(-1,1)\times R$. Let $p$ be a point on $R$ that is disjoint from the projection of $L$ to $R$. In \cite{li2022instanton}, the authors studied the instanton homology group 
$$\SiHI_{R,p}(L) = \II(S^1\times R, L, S^1\times\{p\}|\{t_*\}\times R),$$ where $S^1$ is viewed as the quotient space of $[-1,1]$ with $-1$ identified with $1$, and $t_*\in S^1$ is a fixed base point. 

Suppose $c$ is an embedded $1$--manifold in $R$, and $b$ is a band attached to $c$ that is disjoint from $p$. Then the band surgery along $b$ defines a link cobordism from $c$ to $c_b$ as links in $(-1,1)\times R$. Therefore, it induces a cobordism map for Floer homology groups (up to sign)
$$
\SiHI(b): \SiHI_{R,p}(c) \to \SiHI_{R,p}(c_b).
$$

As discussed in \cite[Proposition 6.12]{li2022instanton}, the maps $\SiHI(b)$ are components of the second page of a variant of Kronheimer--Mrowka's spectral sequence. 
In \cite[Proposition 6.11]{li2022instanton}, the cobordism maps $\SiHI(b)$ were computed for multiple special cases, and they all have the same structure as $T_\lambda(b)$ up to multiplications by integer powers of $i$ and suitable changes of variables. This computation motivated the definition of the homology invariant $\SiKh_\lambda$.
It is natural to conjecture that the second page of Kronheimer--Mrowka's spectral sequence is isomorphic to the homology $\SiKh_\lambda$ with coefficient ring $\bC$ for some $\lambda\neq 0$. 

In order to establish Theorem \ref{thm_detection}, 
we prove the following technical results that sharpens some of the computations in \cite{li2022instanton}. Let $\lambda_1,\dots,\lambda_4$ be the constants from \cite[Section 6]{li2022instanton}. By  \cite[Lemma 6.9]{li2022instanton}, one can always rescale the generators of the instanton homology groups so that $\lambda_1=\pm1$, $\lambda_3=\pm1$.

\begin{lem}
	Assume the generator $w_0$ defined in \cite[Section 5.2.1]{li2022instanton} is chosen so that $\lambda_1=\pm1$, $\lambda_3=\pm1$. Then $\lambda_2= \pm \lambda_4$. 
\end{lem}

\begin{figure}
	\begin{overpic}[width=0.4\textwidth]{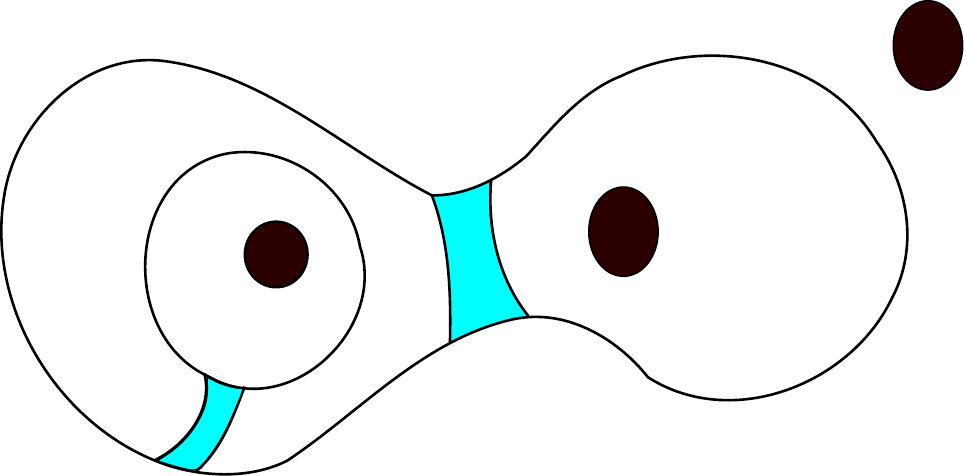}
	\end{overpic}
	\caption{Two bands}\label{fig:TQFT_two_holes}
\end{figure}

\begin{proof}
	Consider the two bands in Figure \ref{fig:TQFT_two_holes} and apply the TQFT property of $\SiHI(b)$. 
\end{proof}

\begin{figure}
	\begin{overpic}[width=0.4\textwidth]{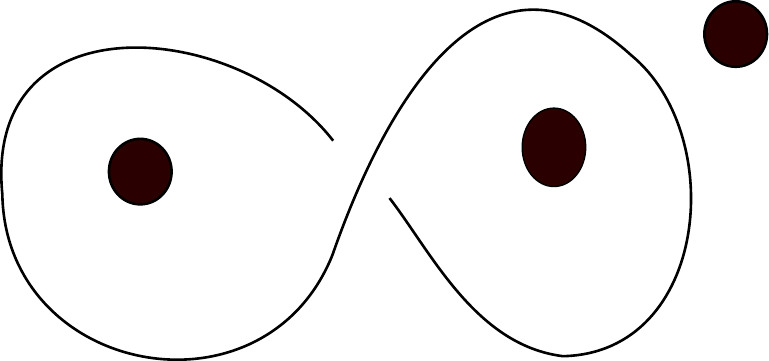}
	\end{overpic}
	\caption{The knot $K$ in $(-1,1)\times \Sigma$}\label{fig:eight_knot_Sigma}
\end{figure}

\begin{figure}
	\begin{overpic}[width=0.4\textwidth]{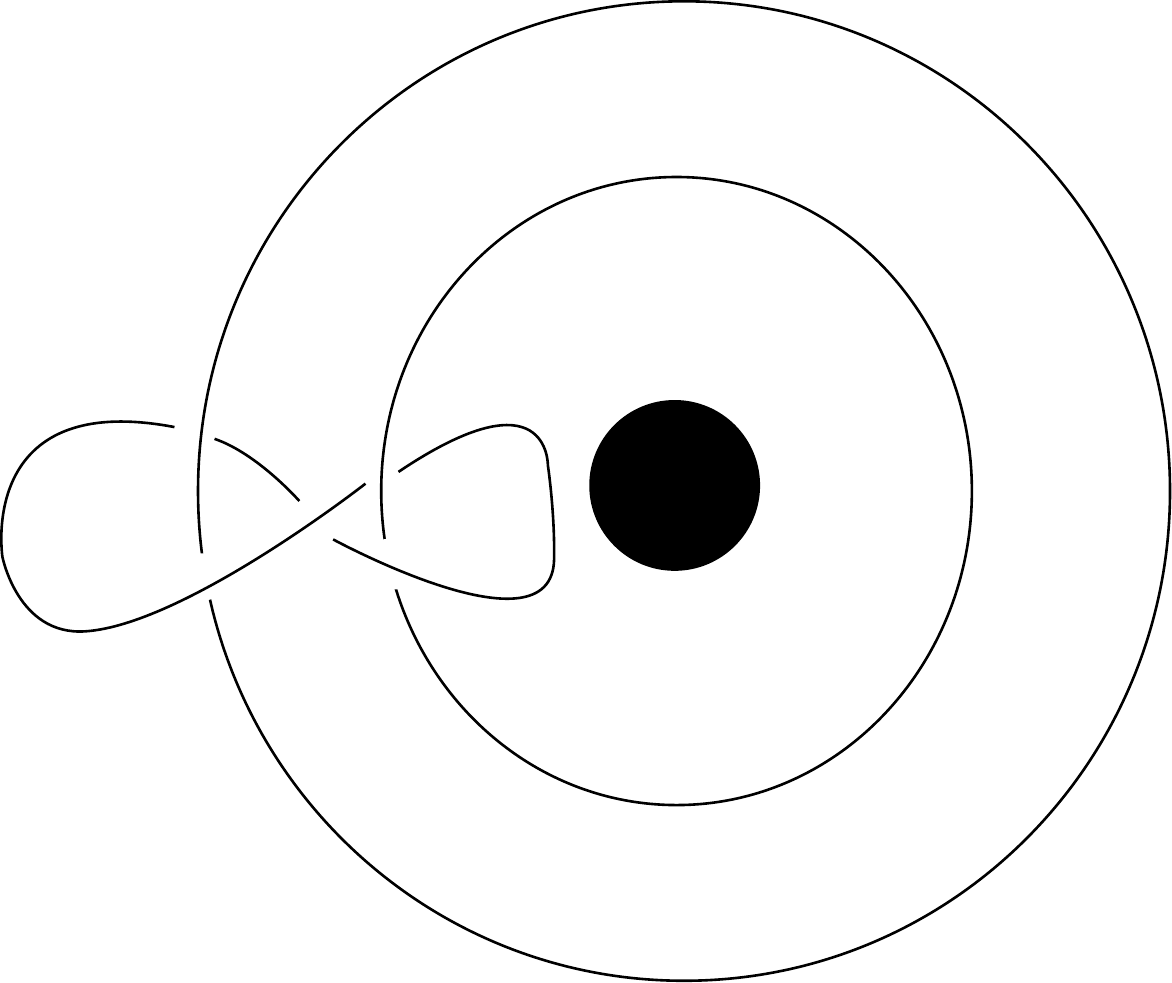}
	\end{overpic}
	\caption{The annular link $L$}\label{fig:eight_annular}
\end{figure}

\begin{lem}
	The coefficients $\lambda_2$ and $\lambda_4$ are both non-zero.
\end{lem}

\begin{proof}
	We keep using the notation from \cite[Section 6]{li2022instanton}. According to the proof of \cite{li2022instanton}*{Lemma 6.4}, we see that it suffices to prove that  
	$$
	\SHI([-1.1]\times \Sigma, \{0\}\times \Sigma, K)\le 4.
	$$
	where $K$ is the knot shown in Figure \ref{fig:eight_knot_Sigma}. 
	
	Let $L$ be the link in the thickened annulus as shown in Figure \ref{fig:eight_annular}.
	Pick a meridional disk in the thickened annulus which intersects $L$ at two points. We decompose the thickened annulus
	along this disk and obtain a product sutured thickened disk with a tangle $T$ in it. The sutured intanton Floer homology
	of this sutured manifold with tangle $T$ is isomorphic to $\AHI(L,2)$ according to \cite{LXZ-unknot}*{Theorem 2.14}, where $\AHI(L,2)$ denotes the component of the annular instanton Floer homology with Alexander grading $2$.
	The tangle 
	$T$ has two product vertical components.
	We remove the tubular neighborhoods of the 
	two vertical components and add a meridian suture to the boundary of each neighborhood to obtain a sutured manifold
	$M'$ with a knot $K'$ in it. Moreover, this process does not change the sutured instanton Floer homology
	according to \cite{XZ:excision}*{Lemma 7.10} and its proof. Therefore we have
	$$
	\SHI(M',\gamma_{M'}, K')\cong \AHI(K,2).
	$$
	Notice that in the definition of sutured instanton Floer homology, the pair $(M',K')$ and $(M,K)$ can be given the same closure, therefore their sutured instanton homology are isomorphic. As a result, we have
	\begin{equation}\label{eq_SHI=AHI}
		\SHI([-1,1]\times \Sigma, \{0\}\times \Sigma, K)\cong \AHI(L,2)
	\end{equation}
	A straightforward calculation shows that 
	$$
	\AKh(L,2;\mathbb{C})\cong\mathbb{C}^4
	$$
	where $\AKh(L,2;\mathbb{C})$ denotes the component of the annular Khovanov homology of $L$ with Alexander grading $2$ and with coefficient ring $\bC$.
	According to \cite{AHI}*{Theorem 5.16}, we have
	$$
	\dim \AHI(L,2;\mathbb{C})\le \dim \AKh(L,2;\mathbb{C})=4.
	$$
	Therefore Equation \eqref{eq_SHI=AHI} implies 
	$$
	\dim \SHI([-1,1]\times \Sigma, \{0\}\times \Sigma, K)\le 4.
	$$ 
	So we obtain $\lambda_2\neq 0$. Since $\lambda_4=\pm \lambda_2$, we also have
	$\lambda_4\neq 0$.
\end{proof}

\begin{proof}[Proof of Theorem \ref{thm_detection}]
	Recall that when $\Sigma$ is a compact surface with genus zero, there is a grading on $V(-)$ such that $T_0(b)$ is homogeneous with degree zero and $T_1(b)$ is homogeneous with degree $-1$. Since $\lambda_2\neq 0$, we can rescale the map $\Theta_{w_0, \sigma}$ in \cite{li2022instanton} by a factor of $\lambda_2^k$ at degree $k$. By the discussions in \cite{li2022instanton}*{Section 6}, there is a spectral sequence of chain complexes in $\bC$--coefficients that converges to $\II([-1,1]\times \Sigma, \{0\}\times \partial \Sigma, L)$, whose second page $(E_2,d_2)$ is isomorphic to the chain complex $(\CKh_{\Sigma,1}(L), \mathcal{D}_{\Sigma,1})$ up to multiplications by integer powers of $i$ on the components of the differential map. In other words, there exists a chain complex $(C,d)$ defined with $\bZ[i]$ coefficients, such that when reducing to $\bC$ coefficients, it is isomorphic to $(E_2,d_2)$; when reducing to $\bZ[i]/(i-1)\cong \bZ/2$ coefficients, it is isomorphic to the chain complex $(\CKh_{\Sigma,1}(L), \mathcal{D}_{\Sigma,1})$.
	By the universal coefficient theorem, we have  
	\begin{align*}
	\rank_{\bZ/2} \SiKh_{\Sigma,1}(L;\bZ/2)
	&\ge \rank_{\bZ[i]}H(C,d) 
	= \dim_\bC H(E_2,d_2) 
	\\
	&\ge \dim_\bC\II([-1,1]\times \Sigma, \{0\}\times \partial \Sigma, L),
	\end{align*}
	and the desired result follows from \cite[Theorem 1.3]{li2022instanton}.
\end{proof}

\bibliographystyle{amsalpha}
\bibliography{references}

\end{document}